 \documentclass[final,3p,times,12pt]{elsarticle}

\usepackage{amssymb}
\usepackage{amsmath}
\usepackage{hyperref}
\usepackage{pstricks,pst-plot,pst-node} 
\usepackage{infix-RPN}

\journal{Nuclear Physics B}

\newtheorem{theorem}{Theorem}

\newtheorem{proposition}{Proposition}

\newtheorem{definition}{Definition}
\newtheorem{remark}{Remark}
\newtheorem{algo}{Algorithm}

\newcommand{\half}{\tfrac{1}{2}}

\newcommand{\ext}{\text{ext}}
\newcommand{\Hd}{{\mathbf H}_\text{d}}
\newcommand{\Ha}{{\mathbf H}_\text{aniso}}
\newcommand{\Heff}{{\mathbf H}_\text{eff}}
\newcommand{\bH}{\bar {\mathbf H}_\text{eff}}
\newcommand{\He}{{\mathbf H}_\text{ext}}
\newcommand{\tf}{{\mathbf \Psi}}
\renewcommand{\d}{\partial}

\newcommand{\dx}{\,{\operatorname d} {\mathbf x}}
\newcommand{\dt}{\,{\operatorname d} t}
\renewcommand{\d}{\,{\operatorname d}}
\newcommand{\Om}{\Omega}
\newcommand{\OmT}{{\Omega_T}}
\newcommand{\OT}{{\Omega_T}}

\newcommand{\vecm}{{\mathbf m}}
\newcommand{\vece}{{\mathbf e}}
\newcommand{\vecv}{{\mathbf v}}
\newcommand{\vecu}{{\mathbf u}}
\newcommand{\vecx}{{\mathbf x}}
\newcommand{\vecw}{{\mathbf w}}

\newcommand{\Div}{\text{div }}

\begin{document}
 
\begin{frontmatter}



\author[cmap]{Fran\c{c}ois Alouges\corref{cor1}}
\ead{francois.alouges@polytechnique.edu}
\author[inp1,inp2]{Evaggelos Kritsikis}
\author[cmap]{Jutta Steiner}
\author[inp1,inp2]{Jean-Christophe Toussaint}
\cortext[cor1]{Corresponding author}
\address[cmap]{CMAP - Ecole Polytechnique, Route de Saclay 91128 Palaiseau cedex, France}

\title{A convergent and precise finite element scheme for Landau-Lifschitz-Gilbert equation}



\address[inp1]{Institut N\'eel, 25 avenue des Martyrs, B\^atiment K, BP 166, 38042 Grenoble cedex 9, France}
\address[inp2]{Grenoble-INP, 46 Avenue FŽlix Viallet, 38031 Grenoble Cedex 1, France}

\begin{abstract}
  In this paper, we rigorously study an order 2 scheme that was
  previously proposed by some of the authors. A slight modification is
  proposed that enables us to prove the convergence of the scheme
  while simplifying in the same time the inner iteration.

\end{abstract}


\end{frontmatter}






\section{Introduction}
In 1935 Landau and Lifschitz proposed an equation that models the
magnetization in a ferromagnetic material
\cite{LandauLifschitz}. Supposing that the three dimensional
ferromagnetic sample occupies some domain $\Omega\subset \mathbb{R}^3$
and calling $\vecm$ the {\it direction} of the magnetization, the
Landau-Lifschitz-Gilbert (LLG) equation reads
\begin{equation}
\left[
\begin{array}{l}
  \displaystyle \partial_t \vecm -\alpha\, \vecm\times \partial_t \vecm=- \gamma_0\vecm\times \Heff\text{ in }\Omega ,\\
  \displaystyle \partial_n \vecm=0\text{ on }\partial \Omega.
\end{array}
\right.
\label{LL1}
\end{equation}
The parameters in the equation are the damping parameter $\alpha$ and
the gyromagnetic constant $\gamma_0$. The so-called effective magnetic
field $\Heff$ is given by the functional derivative of the
micromagnetic (free) energy $\mathcal{E}$, more precisely
\begin{equation}
\Heff(\vecm)=-\frac{\partial \mathcal{E}}{\partial \vecm}
=d^2\,\Delta \vecm+\Hd(\vecm)
+ \He+Q\,(\vece\cdot \vecm)\,\vece
\label{heff}
\end{equation}
where the energy $\mathcal{E}$ (see
\cite{LandauLifschitz,Brown,HubertSchafer}) is given by
\begin{equation}
\mathcal{E}(\vecm)=\frac12\left(\,d^2\,\int_\Omega |\nabla \vecm|^2\d\vecx 
-\int_\Omega \Hd(\vecm)\cdot \vecm\d\vecx 
-2 \int_\Omega \He\cdot \vecm\d\vecx 
- Q\int_\Omega (\vece\cdot \vecm)^2\d\vecx\right).
\label{E(m)}
\end{equation}
The four contributions to the effective field in \eqref{heff} and the
energy in \eqref{E(m)}, respectively, correspond to the so-called
exchange, stray-field, applied and anisotropy field or energy,
respectively. The material constants in \eqref{heff} and \eqref{E(m)}
are the exchange constant $d$, the anisotropy constant $Q$ and the
anisotropy direction $\vece$ (also called the easy axis).
Furthermore, the vector field $\He$ models an applied magnetic
field. We will also use the notation $\Ha=Q\,(\vece\cdot \vecm)\,
\vece$.  The stray field $\Hd(\vecm)$ is the magnetic field induced by
the magnetization distribution $\vecm$ via the following (subset of)
static Maxwell equations
\begin{equation}
\left[
\begin{array}{l}
\text{curl }\Hd (\vecm)=0\text{ in }\mathbb{R}^3\\
\Div (\Hd (\vecm)+\vecm)=0\text{ in }\mathbb{R}^3\,.
\end{array}
\right.
\label{Maxwell}
\end{equation}
Below the Curie temperature, the magnetization can be described by a
directional field that we rescale to be of unit length. It is
straightforward to check that the magnitude of the magnetization
\begin{equation}
|\vecm(\vecx,t)|=1
\label{constraint}
\end{equation}
is conserved by the dynamics \eqref{LL1}. Take note that the
gyromagnetic term is a conservative term while the damping term leads
to the following energy dissipation law
\begin{equation}
  \frac{\d}{\dt}\mathcal{E}(\vecm(t))=-\frac{\alpha}{\gamma_0} \int_\Omega |\partial_t \vecm|^2\d\vecx\,.
\label{dissipation}
\end{equation}
Rescaling time and redefining $\alpha$ allows to assume that
$\gamma_0=1$.  \medskip

The numerical approximation of solutions to \eqref{LL1} is an
important issue in applications. Nowadays, numerous strategies exist
in the literature -- among them only few reliable ones. Classical
schemes are based on finite differences that, as usual, are well
adapted to Cartesian grids. On the other hand, finite elements
approximations are well suited in case of complex geometries and weak
solutions, though bearing the drawback that they are in practice
difficult to analyze. In particular, proving the convergence of a
finite element solution towards a solution of \eqref{LL1} as the space
and time steps tend to zero turns out to be quite difficult and has
probably been first established in \cite{ALOUGESJAISSON}. This result
was further improved in \cite{BARTELSPROHL} and \cite{alougesAIMS},
for the case where only the exchange term is present. We hereafter
study a further generalization of the scheme proposed in
\cite{alougesAIMS}: An order 2 (in time) variant. Numerical tests
support the performance of the method.
\medskip

Let us start with brief outline of our paper. In Section \ref{sec2} we
first recall the notion of weak solutions. Section \ref{sec3}
introduces the finite elements spaces. Section \ref{revisiting}
restates the order one scheme as proposed in \cite{alougesAIMS}. The
nonlinearity of the LLG equation calls for recurrent renormalization
of the time-discrete approximation. This issue is also discussed in
Section \ref{sec5}.  Section \ref{sec6} finally provides a derivation
of our new scheme, the main result about its convergence and its
proof.
\section{Notion of weak solutions to LLG}
\label{sec2}
Let us recall the notion of a weak solution
to \eqref{LL1} from \cite{ALOUGESSOYEUR} and \cite{VISINTIN}.
\begin{definition} \label{timeinterpol} Consider an initial
  magnetization, i.e., a vector field $\vecm_0 \in H^1(\Omega)^3$ that
  is a.e.\ of unit length.  A vector field $\vecm$ is called a weak
  solution to \eqref{LL1} with initial data $\vecm_0$ if for all times
  $T>0$ there holds
\begin{enumerate}
\item $\vecm\in H^1(\OT)^3$ with $\OT=\Omega\times(0,T)$, and $|\vecm|=1$
  a.e.\
\item for all test functions $\tf\in H^1(\OT)^3$
\begin{multline}
  \int_{\OT}\partial_t \vecm \cdot \tf \,\dx\,\dt - \alpha \int_{\OT}
  \left(\vecm \times \partial_t \vecm \right) \cdot \tf \,\dx\,\dt \\
  = \,d^2\,\sum_{i=1}^d \int_{\OT}\left(\vecm
    \times \partial_{x_i}\vecm \right) \cdot \partial_{x_i}\tf
  \,\dx\,\dt -\int_{\OT}\vecm\times (\Hd(\vecm)+\He+\Ha(\vecm))\cdot
  \tf \,\dx\dt,
  \label{ws}
\end{multline}
\item the magnetization initially satisfies
  $\vecm(\vecx,0)=\vecm_0(\vecx)$ in the trace sense, and
\item the energy decreases according to
\begin{equation}
  \mathcal{E}(\vecm(T))+\alpha\int_{\OT}\left|\partial_t \vecm\right|^2\dx\,\dt 
\le \mathcal{E}(\vecm(0)).
\label{energyestimate}
\end{equation}
\end{enumerate}
\end{definition}
\section{The finite element scheme}
\label{sec3}
As in \cite{ALOUGESJAISSON}, our discretization relies on piecewise
linear finite elements in space combined with a linear interpolation
in time. The domain $\Omega$ is discretized by a conformal
triangulation $\mathcal{T}_h$ of mesh size $h$ with vertices
$(\vecx_i^h)_{1\leq i\leq N_h}$. Let us denote by $(\phi_i^h)_{1\leq i\leq
  N_h}$ the set of associated piecewise linear basis functions that
satisfy $\phi_i^h(\vecx_j^h)=\delta_{i,j}$ at the vertices $\vecx_j^h$ for
${1\leq i,j \leq N_h}$, where $\delta_{i,j}$ denotes the Kronecker
symbol. This amounts to a standard
$P^1(\mathcal{T}_h)$-discretization. Based on the scalar basis
$(\phi_i^h)_{1\leq i\leq N_h}$ we construct the {\it vector}-valued
finite element space in the form of
\begin{equation*}
  V_h=\left\{{\mathbf u}^h=\sum_i{\mathbf u}_i \phi_i^h,\,\text{ s.t.\ }\forall i,\,{\mathbf u}_i\in \mathbb{R}^3\right\}.
\end{equation*}
Due to the constraint \eqref{constraint}, the solution to \eqref{ws} is sought for in the subset
\begin{equation*}
  M_h=\left\{{\mathbf u}^h\in V_h,\,\text{ s.t.\ }\forall i,\,{\mathbf u}_i\in \mathbb{S}^2\right\} \subset V_h.
\end{equation*}
Let us also introduce the tangent space in $\vecm^h=\sum_i \vecm_i \phi_i^h\in M_h$ is denoted by
\begin{equation*}
K_\vecm=\left\{{\mathbf v}^h=\sum_i {\mathbf v}_i \phi_i^h,\,\text{ s.t.\ }\forall i,\,{\mathbf v}_i\cdot \vecm_i=0\right\}.
\end{equation*}
Furthermore, the classical nodal interpolation operator is given by
\begin{align}
  \mathcal{I}_h:{\mathcal C}^0(\Omega,\mathbb{R}^3)&\rightarrow V_h\notag\\
 {\mathbf u}&\mapsto \sum_i{\mathbf u}(\vecx_i^h)\phi_i^h.\label{JS:interpol}
\end{align} 
To simplify notations, the index $h$ of the ansatz functions will be
neglected from now on most of the times, i.e., we write ${\mathbf u}$,
${\mathbf v}$, etc.\ instead of ${\mathbf u}^h$, ${\mathbf v}^h$,
respectively, in case this does not lead to any ambiguities.
\section{Revisiting the $\theta$-scheme}
\label{revisiting}
The finite element scheme proposed in \cite{ALOUGESJAISSON} relies on
the observation that the LLG equation \eqref{LL1} -- with the notation
$\displaystyle \vecv=\partial_t \vecm$ -- can be rewritten in the following weak form
\begin{multline}
  \alpha\int_\Omega \vecv\cdot \tf\dx+ \alpha\int_\Omega \vecm\times \vecv \cdot \tf\dx
 \\ =-\,d^2\,\int_\Omega\nabla \vecm\cdot \nabla \tf \dx+\int_\Omega(\Hd(\vecm)+\He+\Ha(\vecm)) \cdot \tf\dx.
\label{disLLG}
\end{multline}
Equation \eqref{disLLG} holds for every test function $\tf\in H^1(\Omega,\mathbb{R}^3)$
that satisfies $\tf(\vecx)\cdot \vecm(\vecx)=0$ for a.e. $\vecx$ in
$\Omega$.  The reformulation of \eqref{LL1} in the form of
\eqref{disLLG} motivated the following $\theta-$ scheme introduced in
\cite{alougesAIMS}:
\begin{algo}
  \label{algo1}
  Given an initial $\vecm^0\in M_h$ choose $\theta\in[0,1]$ and a time
  step size $\tau=\frac{T}{N}$ with $N\in \mathbb N$. For
  $n=0,1,\dots,N$
 \begin{equation}
\begin{array}{l} 
 \left[
\begin{array}{l}
  \text{a) find } \vecv^n\in K_{\vecm^n}\text{ such that for all test functions } \tf\in K_{\vecm^n}\\\\
  \displaystyle \alpha \int_\Omega \vecv^n\cdot \tf\dx+\int_\Omega \vecm^n\times \vecv^n\cdot \tf \dx
  \\\displaystyle\quad=-d^2\int_\Omega\,\nabla (\vecm^n+\theta \tau \vecv^n)\cdot \nabla \tf\dx
  +\int_\Omega  (\Hd(\vecm^n)+\He+\Ha(\vecm^n)) \cdot \tf\dx\,\\\\
  \displaystyle \text{b) set }\vecm^{n+1}=\sum_i \vecm^{n+1}_i\phi_i^h,\text{ where } \forall i,\,\,\vecm^{n+1}_i=\frac{\vecm^n_i+\tau \vecv^n_i}{|\vecm^n_i+\tau \vecv^n_i|}, 
\end{array}
\right.
\end{array}
\label{theta}
\end{equation}
\end{algo}
It is noteworthy that this procedure requires the solution of a linear
equation in each time step only. Moreover, due to the fact that the
symmetric part of the underlying matrix is positive definite,
existence and uniqueness of a solution to \eqref{theta} is guaranteed.
\medskip

The time discrete solution constructed via algorithm \eqref{theta} at
time-steps $\displaystyle N=\left[\frac{T}{\tau}\right]$ is
interpolated as follows:
\begin{definition}
  In each time interval $t\in [n\tau,(n+1)\tau)$ with $n \in \{0,\cdots,N\}$
  we set
\begin{align*}
  &\vecm_{h,\tau}=\frac{t-n\tau}{\tau}\vecm^{n+1}+\frac{(n+1)\tau-t}{\tau}\vecm^n,\\
  &\vecm_{h,\tau}^-=\vecm^n, \,\,\vecv_{h,\tau}=\vecv^{n}.
\end{align*}
\end{definition}
Our notational convention is thus that $\vecm_{h,\tau}$,
$\vecm_{h,\tau}^-$ and $\vecv_{h,\tau}$ refer to suitable time
interpolants of the time discrete approximation $\vecm^{n}$ and
$\vecv^{n}$.  Notice that $\vecm_{h,\tau}$ is piecewise linear in time
whereas $\vecm_{h,\tau}^-$ and $\vecv_{h,\tau}$ are piecewise
constant. (The introduction of the piecewise constant magnetization
will be useful in the convergence proof.) Based on this
discretization, weak convergence of the constructed approximation was
established in \cite{alougesAIMS}.  Both the proof of this result and
the proof in case of our new scheme consist of the following two main
``classical'' steps: As a first step establishing an energy estimate
which guarantees the convergence (sufficiently strong) of the sequence
constructed and then in a second step verifying that the limit indeed
satisfies the equation. As far as the first step is concerned, the
following section addresses the fact that the energy behaves well
under renormalization -- in principle a strongly nonlinear modification of the flow.
\section{Renormalization decreases the energy}
\label{sec5}
The influence of the renormalization on the {\it exchange}
energy was for instance investigated in \cite{ALOUGES} in the
continuous case. More precisely, it was shown that for maps $\vecw\in
H^1(\Omega,\mathbb{R}^3)$ with $|\vecw(\vecx)|\geq 1$ a.e. $\vecx \in \Omega$ one
has
\begin{equation}
  \int_\Omega \left|\nabla \frac{\vecw}{|\vecw|}\right|^2\dx \leq \int_\Omega \left|\nabla \vecw\right|^2\dx.
\label{energdecay}
\end{equation}
Hence, the renormalization step is expected to be energy decreasing --
a least as far as the Dirichlet energy is concerned. Applications more
related to finite element approximation of micromagnetic
configurations can be found in \cite{ALOUGESETAL}. The discrete
version of \eqref{energdecay} was proved by Bartels in \cite{BARTELS}:
\begin{theorem}
  \label{theobart}{\cite{BARTELS}}
  If the basis functions of the $P^1$-approximation satisfy
  \begin{equation}
    \forall i\ne j,\,\,\int_\Omega\nabla \phi_i^h \cdot \nabla \phi_j^h \dx \leq 0,
    \label{condsuff}
  \end{equation}
  then for all $\vecv=\sum_i \vecv_i \phi_i^h \in V_h$ such that $\forall
  i\in\{1,\cdots,N_h\},\,|\vecv_i|\geq 1$ it holds that
  \begin{equation}
  \int_\Omega \left|\nabla \mathcal{I}_h\left( \frac{\vecv}{|\vecv|}\right)\right|^2 \dx\leq \int_\Omega |\nabla \vecv|^2\dx.
  \label{condtriang}
  \end{equation}
\end{theorem}

In 3d
, the condition \eqref{condsuff} -- and hence \eqref{condtriang}
-- is for instance satisfied provided all dihedral angles of the
tetrahedra of the mesh are smaller than $\pi/2$, see \cite{VANSELOW}.

\section{The new (almost) order $2$-scheme}
\label{sec6}
Let us embark on the motivation and description of our new scheme.  As
remarked in \cite{alougesAIMS}, it is {\it not} sufficient to choose
$\theta=\frac12$ in \eqref{theta} to achieve quadratic order due to the
renormalization which inherently introduces an error of order
2. Hence, it is necessary to modify the time-discrete approximation of
the magnetization $\vecm$. 
\medskip

Consider an iterate $\vecm(n\tau)$ at time
$n\tau$. It is well known that the mid-point rule is exact up to cubic error, i.e., 
\begin{align*}
  \vecm((n+1)\,\tau)=\vecm(n\tau)+\tau
  \vecm_t((n+\tfrac{1}{2})\,\tau)+O(\tau^3).
\end{align*}
Now, given a current iterate $\vecm(n\tau)$ at time $n\tau$, a Taylor
expansion up to cubic order, i.e.,
\begin{equation*}
  \vecm((n+1)\,\tau)=
\vecm(n\tau)+\tau \vecm_t(n\tau)+\frac{\tau^2}{2} \vecm_{tt}(n\tau)+O(\tau^3)
\end{equation*}
reveals that the parallel component of the subsequent iterate (along
$\vecm(n\tau)$) is due to the unit length constraint given by
\begin{equation*}
  \vecm\big((n+1)\,\tau\big)\cdot \vecm(n\tau) = 1 - \tau^2| \vecm_t(n\tau)|^2+O(\tau^3).
\end{equation*}
This can easily be inferred from the unit length constraint by
differentiation, i.e., using the relations
\begin{align*}
  \vecm \cdot \vecm_t=0,\\
  \vecm_t \cdot \vecm_t+\vecm\cdot \vecm_{tt}=0.
\end{align*}
We therefore propose to modify the original first order scheme by
replacing the tangential update with the following higher order
approximation
\begin{align}
  \vecv&=P_{\vecm^\perp}\, \vecm_t((n+\tfrac{1}{2})\,\tau)\notag\\
  &=P_{\vecm^\perp} \,(\vecm_t(n\tau) + \tfrac{\tau}{2}\vecm_{tt}(n\tau)) +O(\tau^2)\notag\\
  &= \vecm_t(n\tau) + \tfrac{\tau}{2}
  \,P_{\vecm^\perp}\vecm_{tt}(n\tau) +O(\tau^2),\label{defv}
\end{align}
where $P_{\vecm^\perp}$ denotes the projection onto the orthogonal
component of $\vecm(n\tau)$. 
\medskip

We will use the short hand notation $\vecm=\vecm(n\tau)$ and
$\vecm_t=\vecm_t(n\tau)$ -- provided that what is stated remains clear
without ambiguity. Let us proceed with the derivation of the equation
that is satisfied by $\vecv= \vecm_t(n\tau) + \tfrac{\tau}{2}
\,P_{\vecm^\perp}\vecm_{tt}(n\tau)$, i.e. the counterpart to
\eqref{disLLG}. The equation will be inferred from the differentiated
LLG equation which we restate as
\begin{align}
  \label{JS:LLG1}
  \alpha \vecm_t+\vecm\times \vecm_t=\Heff(\vecm)-(\Heff(\vecm)\cdot
  \vecm) \,\vecm
\end{align}
by multiplying \eqref{LL1} with $\vecm \times$.  
\begin{remark}
\label{remark1}
Although the mid-point rule is of order $2$, our scheme will be only
almost of order $2$ -- as the section's title suggests and as we will
see in the sequel. We have to introduce a regularizing term in order
to obtain the necessary estimates in the convergence proof. This term
prevents the scheme from being of order 2, in the sense that the
consistency error is not of order $O(\tau^3)$ but only
$O(\tau^{3-\epsilon})$ for any $\epsilon >0$. On the other hand, this
regularization approach allows for {\it unconditional convergence} of
the scheme. If we do not insist on {\it unconditional convergence},
then under the condition $\tau\ll h$, consistency up to order
$O(\tau^3)$ is attainable.
\end{remark}
\medskip

To begin with, the differentiation of \eqref{JS:LLG1} w.r.t. time
yields
\begin{align} 
  \alpha \vecm_{tt}+\vecm_t&\times \vecm_{tt}\\&=\frac{\partial
    \Heff}{\partial \vecm}(\vecm_t)-\left(\frac{\partial
      \Heff}{\partial \vecm}(\vecm_t)\cdot \vecm\right)\, \vecm-
  (\Heff(\vecm)\cdot \vecm_t)\, \vecm-(\Heff(\vecm)\cdot \vecm)\,
  \vecm_t, \label{JS:LLG2}
\end{align}
where 
\begin{equation*}
  \frac{\partial \Heff}{\partial \vecm}
  =d^2\,\Delta \vecm_t+\Hd(\vecm_t)+Q\,(\vece\cdot \vecm_t)\,\vece
\end{equation*}
and where we once again used the unit length constraint
\eqref{constraint}. The application of the projection to
\eqref{JS:LLG2} in combination with \eqref{JS:LLG1} yields
\begin{multline*}
  \int_\Om \alpha   \vecv\cdot \tf +\vecm\times \vecv\cdot \tf\dx\\=\int_\Om \Heff(\vecm)\cdot \tf \dx +\tfrac{\tau}{2} \int_\Om \frac{\partial \Heff}{\partial \vecm}(\vecm_t)\cdot \tf\dx-\tfrac{\tau}{2} \int_\Om( \Heff(\vecm)\cdot \vecm)\,\vecm_t\cdot\tf\dx
\end{multline*}
for any test function $\tf$ with $\tf\cdot \vecm=0$. Observe that
$\vecm_t(n\tau)=\vecv+O(\tau)$, cf. \eqref{defv}. Therefore up to
higher order terms
\begin{multline}
  \int_\Om  (\alpha +\tfrac{\tau}{2} ( \Heff(\vecm)\cdot \vecm) )\,\vecv\cdot \tf+\vecm\times \vecv\cdot \tf\dx -\tfrac{\tau}{2} \int_\Om \frac{\partial \Heff}{\partial \vecm}(\vecv)\cdot \tf\dx\\=\int_\Om \Heff(\vecm)\cdot \tf \dx+O(\tau^2),
\label{JS:inter1}
\end{multline}
where we remind that $\vecm=\vecm(n\tau)$ and
$\vecm_t=\vecm_t(n\tau)$. Observe that the latter equation is (at
first sight surprisingly) {\it linear} in $\vecv$. However, nothing
can be stated about its well-posedness since both the first and the
last contribution on the l.h.s.\ of \eqref{JS:inter1} potentially
affect the definiteness of the symmetric part of the operator. In
order to guarantee solvability and uniqueness we proceed with higher
order modifications that will finally lead to a well posed
formulation. We address the first contribution and define
\begin{align}
\label{JS:inter2}
 \tilde\varphi_M(x)=
  \begin{cases}
    \alpha +\tfrac{\tau}{2} \min(x,M)\text{ for } x\geq0,\\
    \displaystyle{\frac{\alpha}{1+\tfrac{\tau}{2} \min(-x,M)}}\text{
      for } x<0.
  \end{cases}
\end{align}
Notice that $\tilde\varphi_M(x)=\alpha +\tfrac{\tau}{2}
\min(x,M)+O(\tau^2M^2)$. By abuse of notation we define
\begin{equation}
\label{JS:inter3}
\varphi_M(\vecm)=\tilde \varphi_M(\Heff(\vecm)\cdot \vecm).
\end{equation}
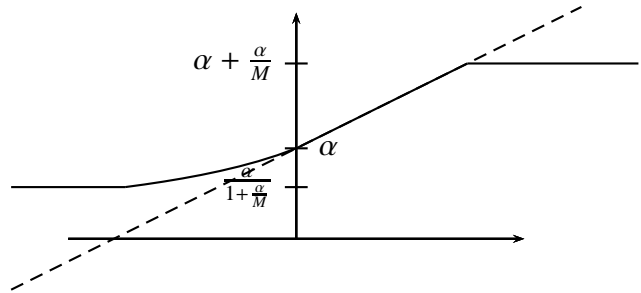
\begin{figure}[htb]
  \begin{center}
\psset{unit=1.5}
\begin{pspicture}(-2,-.5)(2,2) \psset{plotpoints=10}
  \psaxes[labels=none,ticks=none]{->}(0,0)(-2,0)(2,2)
  \psline(-.1,.8)(.1,.8) \psline(-.1,1.55)(.1,1.55)
  \psline(-.1,0.457)(.1,0.457)
  \infixtoRPN{.8+.5*x}\psplot{0}{1.5}{\RPN}
  \infixtoRPN{.8/(1-.5*x)}\psplot{-1.5}{0}{\RPN}
  \infixtoRPN{.8+.5*x}\psplot[linestyle=dashed]{-2.5}{2.5}{\RPN}
  \infixtoRPN{.8+.5*1.5}\psplot{1.5}{3}{\RPN}
  \infixtoRPN{.8/(1+.5*1.5)}\psplot{-1.5}{-2.5}{\RPN}
  \rput[r](-0.2,1.55){$\alpha+\frac{\alpha}{M}$}
  \rput[l](0.2,.8){$\alpha$}
  \rput[r](-0.2,.457){$\frac{\alpha}{1+\frac{\alpha}{M}}$}
\end{pspicture} 
   \end{center}
\caption{The regularizing cut-off function $\tilde \varphi(x)$.}
\label{}
\end{figure}

As long as $\Heff(\vecm)\cdot \vecm $ is uniformly bounded, we derive
from \eqref{JS:inter1} by plugging in \eqref{JS:inter3} that
\begin{equation}
  \int_\Om  \varphi_M( \vecm) \,\vecv\cdot \tf+\vecm\times \vecv\cdot \tf\dx
  -\tfrac{\tau}{2} \int_\Om \frac{\partial \Heff}{\partial \vecm}(\vecv)\cdot
  \tf\dx=\int_\Om \Heff(\vecm)\cdot \tf \dx+O(\tau^2).
\label{JS:inter4}
\end{equation} 
Replacing $\Heff$ and
$\frac{\partial \Heff}{\partial \vecm}$ by their very definition, we obtain
the counterpart to \eqref{disLLG} for our new second order scheme: 
  \begin{multline}
    \int_\Om \varphi_M( \vecm) \,\vecv\cdot \tf+\vecm\times \vecv\cdot
    \tf\dx +\tfrac{\tau}{2} \int_\Om \,d^2\,\nabla \vecv \cdot
    \nabla\tf -\Hd(\vecv)\cdot \tf
    -Q(\vece\cdot \vecv)(\vece \cdot \tf) \dx\\
    =\int_\Om -\,d^2\,\nabla \vecm \cdot \nabla\tf +\Hd(\vecm)\cdot
    \tf +Q(\vece\cdot \vecm)(\vece \cdot \tf) +\He \cdot \tf\dx.
\label{JS:inter5}
\end{multline}
We introduce only one further, final modification which implements the
strategy delineated in Remark \ref{remark1}: In order to maintain
unconditional convergence we additionally modify the second highest order term
on the r.h.s. in the following way
\begin{equation*}
  \tfrac{\tau}{2} \int_\Om \,d^2\,\nabla \vecv \cdot \nabla\tf  \d x 
  \quad\leadsto\quad 
  \tfrac{\tau}{2} \int_\Om (1+\rho(\tau)) \,d^2\,\nabla \vecv \cdot \nabla\tf \dx, 
\end{equation*}
where $\rho(\tau)\rightarrow 0$ as $\tau\rightarrow 0$. Take note
that for $\rho$ decreasing at least linearly, quadratic order is
conserved. However, only in case that $\rho$ is slightly sublinear,
for example $\rho(\tau)=\tau|\ln(\tau)|$, do we in fact achieve
unconditional convergence.  \medskip

Adopting Algorithm \ref{algo1}, we arrive at the following scheme:
\begin{algo}
\label{algo2}
Given an initial $\vecm^0\in M_h$ choose a time step size
$\tau=\frac{T}{N}$ with $N\in \mathbb N$ and appropriate $\rho(\tau)$
and $M$, cf.\ Theorem \ref{thetheorem}. For $n=0,1,\dots, N$
 \begin{equation}
\begin{array}{l} 
 \left[
\begin{array}{l}
  \text{a) find } \vecv^n\in K_{\vecm^n}\text{ such that for all test functions } \tf\in K_{\vecm^n}\\\\
  \displaystyle \int_\Om  \varphi_M( \vecm^n) \,\vecv^n\cdot \tf+\vecm^n\times \vecv^n\cdot \tf\dx
  \\
  \quad\displaystyle+\tfrac{\tau}{2} \int_\Om (1+\rho(\tau)) \,d^2\,\nabla \vecv^n \cdot \nabla\tf -\Hd(\vecv^n)\cdot \tf-Q(\vece\cdot \vecv^n)(\vece\cdot \tf) \dx
  \\\displaystyle=\int_\Om -\,d^2\,\nabla \vecm^n \cdot \nabla\tf
  +\left( \Hd(\vecm^n)+\He+\Ha(\vecm^n)\right)
  \cdot \tf\dx. 
  \\\\
  \displaystyle \text{b) set }\vecm^{n+1}=\sum_i \vecm^{n+1}_i\phi_i^h,\text{ where } \forall i,\,\,\vecm^{n+1}_i=\frac{\vecm^n_i+\tau \vecv^n_i}{|\vecm^n_i+\tau \vecv^n_i|}.
\end{array}
\right.
\end{array}
\label{JS:inter5.2}
\end{equation}
\end{algo}
The {\it appropriate} choice of $\rho$ and $M$ can be inferred from
our convergence result, see Theorem \ref{thetheorem}.  \medskip

Let us sum up: The new scheme replaces the search of $\vecv$ as
solution to \eqref{theta} by the search of $\vecv$ as a solution to
\eqref{JS:inter5.2}. Besides this substitution, the algorithm outlined
in Section \ref{revisiting} remains as before in the sense that the
renormalization and the interpolation w.r.t.\ time are left unchanged.
Since equation \eqref{JS:inter5.2} is linear in $\vecv$, our algorithm
is very favorable in practice. \medskip

Before we state our theorem about the convergence let us explicitly
make a statement about its order.
\begin{proposition}
 \label{prop}
 Consider a smooth (in space and time) solution $\vecm$ to
 \eqref{JS:inter5.2} at time $t+\tau$ and a semi-discrete
 time-approximation to $m$ at time $t+\tau$ on the basis of
 \eqref{JS:inter5.2}. More precisely, given $\vecm$ at time $t=n\tau$ determine
 $\vecv= \vecm_t(n\tau) + \tfrac{\tau}{2}
 \,P_{\vecm^\perp}\vecm_{tt}(n\tau)$ as a solution to the variational
 formulation \eqref{JS:inter5.2} with $\rho(\tau)=0$ and $M(\tau)$
 sufficiently large and set
  \begin{equation*}
    \tilde\vecm(\vecx,t+\tau)=\frac{\vecm(\vecx,t)+\tau \vecv(\vecx,t)}{|\vecm(\vecx,t)+\tau \vecv(\vecx,t)|}\text{ for all } \vecx \in \Omega. 
  \end{equation*}
  Then $\tilde \vecm(t+\tau)$ approximates $\vecm(t+\tau)$ up to cubic
  error in $\tau$.
\end{proposition}

\paragraph{Argument for Proposition \ref{prop}}
    The proof is a direct consequence of the Taylor expansion
    performed in \eqref{defv}.
\begin{remark}
  The smoothness of solutions to \eqref{LL1} has been widely studied
  during the course of the past years. In general, the formation of
  singularities cannot be ruled out and we can usually not assume that
  a solution to the initial value problem will be regular. Our
  statement about the order of the approximation is thus only a first
  little step on the way to a proof of the order of convergence, which
  is way beyond the scope of this paper.
\end{remark}
Let us now turn to the convergence result.
\begin{theorem}\label{thetheorem}
  Let $\vecm_0\in H^1(\Omega,S^2)$. Suppose $\vecm^0 \rightarrow
  \vecm_0$ in $H^1(\Omega)$ as $h\rightarrow 0$. If the regular
  sequence of conformal triangulations $(\mathcal{T}_h)_{h>0}$
  satisfies condition \eqref{condsuff}, then the approximation
  $(\vecm_{h,\tau})$ of the sequence constructed via Algorithm
  \ref{algo2} and interpolated according to Definition
  \ref{timeinterpol} converges (up to the extraction of a subsequence)
  weakly in $H^1(\OT)$ to a weak solution $\vecm$ of \eqref{LL1} as
  $h$ and $\tau$ tend to 0 provided $\rho(\tau)\rightarrow_{\tau
    \rightarrow 0} 0$ and one of the two following conditions hold:
  \begin{itemize}
    \item $\tau^{-1} \rho(\tau)\rightarrow_{(h,\tau)\rightarrow 0}
    \infty$ and $\tau M \rightarrow_{(h,\tau)\rightarrow 0} 0 $ or
\item $\rho\equiv 0$ and $\tau \ll h$ as $(h,\tau)\rightarrow 0$.
\end{itemize}
\end{theorem}

\paragraph{Proof of Theorem \ref{thetheorem}} As stated before, the
proof consists of two main steps: Establishing estimates which
guarantee the existence of a sufficiently strong converging
subsequence, and finally proving that the latter converges indeed to a
solution (which satisfies the energy estimate). We will need the
following classical estimate from elliptic regularity theory, namely
\begin{equation}
  ||\Hd(\vecm)||_{L^p(\Omega)}\leq C||\vecm||_{L^p(\Omega)},
\label{estimcst}
\end{equation}
for all $p \in (1,+\infty)$ and for positive constant $C$ which depend
only on $p$. 
\paragraph{Bounds on the sequence}
As we have already observed, the variational formulation in the iteration of \eqref{JS:inter5} possesses a unique solution $\vecv^n$. We test the equation with $\tf=\vecv^n$ itself to find that
  \begin{multline}
    \int_\Om \varphi_M( \vecm^n) \,|\vecv^n|^2\dx +\tfrac{\tau}{2}
    \int_\Om (1+\rho(\tau))\,d^2\,|\nabla \vecv^n |^2
    -\Hd(\vecv^n)\cdot \vecv^n -Q(\vecv^n\cdot e)^2 \dx\\=\int_\Om
    -\,d^2\,\nabla \vecm^n \cdot \nabla \vecv^n +\Hd(\vecm^n)\cdot
    \vecv^n +Q(\vece\cdot\vecm^n)(\vece \cdot\vecv^n) +\He \cdot
    \vecv^n\dx.
\label{JS:inter6}
\end{multline}
Since we assume that the triangulation $\mathcal{T}_h$ satisfies the angle condition \eqref{condtriang} we have that
\begin{align*}
  \int_\Omega \left|\nabla \vecm^{n+1}\right|^2\dx&\leq\int_\Omega
  \left|\nabla (\vecm^n+\tau \vecv^n)\right|^2\dx\\
  &\leq\int_\Omega \left|\nabla \vecm^n\right|^2\dx+2 \tau \int_\Omega
  \nabla \vecm^n \cdot \nabla \vecv^n \dx+ \tau^2\int_\Omega\left|\nabla
    \vecv^n\right|^2\dx.
\end{align*}
Using \eqref{JS:inter6} we obtain that
\begin{multline}
  \,d^2\,\int_\Omega \left|\nabla \vecm^{n+1}\right|^2\dx\leq\,d^2\,\int_\Omega \left|\nabla \vecm^n\right|^2\dx-2\tau\int_\Omega \varphi_M(\vecm^n)\left|\vecv^n\right|^2\dx+\tau^2\int_\Omega \Hd(\vecv^n)\cdot \vecv^n+Q(\vece \cdot\vecv_n)^2\dx\\
  \qquad+2\tau\int_\Omega \Hd(\vecm^n)\cdot \vecv^n+\Ha(\vecv^n)\cdot
  \vecv^n+\He\cdot \vecv^n\dx
  -\tau^2\rho(\tau)\,d^2\,\int_\Omega\left|\nabla \vecv^n\right|^2\dx.
\label{estimenerg}
\end{multline}
Before we move on, let us just rewrite the latter estimate as 
\begin{align}
  \int_\Omega \left|\nabla \vecm^{n+1}\right|^2\dx&\leq\,d^2\,\int_\Omega
  \left|\nabla \vecm^n\right|^2\dx-2\tau\int_\Omega
  \varphi_M(\vecm^n)\left|\vecv^n\right|^2\dx+\tau^2\int_\Omega
  \frac{\partial\bH}{\partial \vecm}(\vecv^n)\cdot \vecv^n\dx\notag\\
  &\qquad+2\tau\int_\Omega \bH(\vecm^n)\cdot \vecv^n\dx
  -\tau^2\rho(\tau)\,d^2\,\int_\Omega\left|\nabla \vecv^n\right|^2\dx.
\label{estimenerg.2}
\end{align}
We partially neglect the negative contributions on the r.h.s.\ of
\eqref{estimenerg.2} -- those which are quadratic in $\vecv^n$ -- and use
\eqref{estimcst} to obtain
\begin{align}
  &\,d^2\,\int_\Omega \left|\nabla \vecm^{n+1}\right|^2\dx\notag\\&\leq\,d^2\,\int_\Omega \left|\nabla \vecm^n\right|^2\dx-2\tau\int_\Omega \varphi_M(\vecm)\left|\vecv^n\right|^2\dx+2\tau||\bH(\vecm^n)||_{L^2(\Omega)}||\vecv^n||_{L^2(\Omega)}-\tau^2\rho(\tau)\,d^2\,\int_\Omega\left|\nabla \vecv^n\right|^2\dx\notag\\
  &\leq \,d^2\,\int_\Omega \left|\nabla
    \vecm^n\right|^2\dx-2\tau\int_\Omega
  \varphi_M(\vecm)\left|\vecv^n\right|^2\dx+ C\tau
  ||\vecv_n||_{L^2(\Omega)}
  -\tau^2\rho(\tau)\,d^2\,\int_\Omega\left|\nabla \vecv^n\right|^2\dx,
\label{estimenerg2}
\end{align}
where the generic constant $C$ depends on $Q$ and $|\Omega|$. Due to
Young's inequality, we have that $C\tau ||\vecv_n||_{L^2(\Omega)}\leq \tau
\beta||\vecv_n||^2_{L^2(\Omega)}+ \frac{\tau C^2}{4\beta}$ for
$\beta>0$. Using the uniform bound $$\varphi_M(\vecm)\geq
\beta=\frac{\alpha}{1+\tfrac{\tau}{2} M}$$ we find by rewriting
\eqref{estimenerg2} that
\begin{multline}
  \,d^2\,\int_\Omega \left|\nabla \vecm^{n+1}\right|^2\dx+ \beta\tau||\vecv^n||^2_{L^2(\Omega)}
  +\tau^2\rho(\tau)\,d^2\,\int_\Omega\left|\nabla \vecv^n\right|^2\dx\\\leq \,d^2\,\int_\Omega \left|\nabla \vecm^n\right|^2\dx+ \frac{\tau C^2(Q, |\Omega|)}{4\beta}.
\label{estimenerg3}
\end{multline}
Summing up in \eqref{estimenerg3} over the time steps we find that
\begin{multline}
  \,d^2\,\int_\Omega \left|\nabla \vecm^{N}\right|^2\dx
  +\beta\tau\sum_{n=0}^{N-1}\int_\Omega
  \left|\vecv^{n}\right|^2\dx+\tau^2\rho(\tau)\,d^2\,\int_\Omega\left|\nabla
    \vecv^n\right|^2\dx\\\leq C\left(T, \,d^2\,\int_\Omega \left|\nabla
    \vecm_0\right|^2\dx, \beta, Q , \Ha\right)
  \label{estimenerg4}
\end{multline}
\medskip

From now on, most of the arguments follow the same line as in
\cite{alougesAIMS}. It holds that 
\begin{equation}
  \left|\frac{\vecm_i^{n+1}-\vecm_i^n}{\tau}\right|\leq |\vecv^n_i|,\text{ for all }n\leq N,\,\, \text{ and }i\in\{1,\cdots,N_h\}.\nonumber
\end{equation}
Moreover, there exists $c>0$ such that for all $1\leq p < +\infty$ and all $\phi_h\in V_h$ there holds
\begin{equation}
\frac{1}{c}||\phi_h||^p_{L^p(\Omega)}\leq h^d\sum_i |\phi_h(x_i^h)|^p \leq c||\phi_h||^p_{L^p(\Omega)},
\label{ineg}
\end{equation}
which implies
\begin{equation}
  \left\Vert \frac{\vecm^{n+1}-\vecm^n}{\tau}\right\Vert_{L^2} \leq c^2 ||\vecv^n||_{L^2}.\label{ineg2}
\end{equation}
Hence we obtain from the energy estimate \eqref{estimenerg4} using \eqref{ineg2} the following bounds
\begin{align}
  \vecm_{h,\tau}\ \text{ is uniformly bounded in } H^1(\OT),
\label{JS:bounds1}\\ 
 \vecv_{h,\tau}\ \text{ is uniformly bounded in } L^2(\OT). \label{JS:bounds2}
\end{align} 
Due to \eqref{JS:bounds1} and \eqref{JS:bounds2}, there exist $\bar\vecm\in H^1(\OT)$ and $\vecv\in L^2(\OT)$ such that up to the extraction of subsequences
\begin{align}
  \label{19}&\vecm_{h,\tau}\rightharpoonup_{(h,\tau)\rightarrow 0} \bar\vecm\text{ weakly in }H^1(\OT),\\
  \label{20}&\vecm_{h,\tau}\rightarrow_{(h,\tau)\rightarrow 0} \bar\vecm\text{ strongly in }L^2(\OT),\\
  \label{21}&\vecv_{h,\tau}\rightharpoonup_{(h,\tau)\rightarrow 0}
  \vecv\text{ weakly in }L^2(\OT).
\end{align}
In addition, we have from \eqref{estimenerg4} that
\begin{equation*}
\sum_{n=0}^{N-1}\tau^2\rho(\tau)\int_\Omega\left|\nabla
  \vecv^n\right|^2\dx= \tau\rho(\tau)\int_0^{T}\int_\Omega\left|\nabla
  \vecv_{h,\tau}\right|^2 \dx\leq C <+\infty
\end{equation*}
If $\rho$ decreases only {\it sublinearly}, i.e. $\tau^{-1} \rho(\tau)\rightarrow_{\tau\rightarrow 0}+\infty$, we deduce that 
\begin{equation}
  \label{bound1}
  \tau \,||\nabla \vecv||_{L^2(\OT)}\rightarrow_{(h,\tau)\rightarrow 0}0.
\end{equation}
If $\rho$ decreases {\it linearly} or faster we have to resort to the
inverse estimate $||\nabla \vecv||_{L^2(\OT)}\lesssim
\frac{1}{h}||\vecv||_{L^2(\OT)}$ in order that estimate \eqref{bound1}
holds true. In fact, is easily seen that \eqref{bound1} is follows
from the inverse estimate in case of $\tau\ll h$.

\paragraph{Preliminary estimates}
We want to prove that $\bar \vecm$ satisfies \eqref{ws} and follow the
strategy of \cite{alougesAIMS}. To begin with, we restate some further
estimates from \cite{alougesAIMS} and derive some necessary statements
about convergence. Observe that for all $n=0,\cdots,J$ and all $t\in
[n\tau,(n+1)\tau)$
\begin{equation}
  |\vecm_{h,\tau}(\vecx,t)-\vecm^-_{h,\tau}(\vecx,t)|=\left|(t-n\tau)\left(\frac{\vecm^{n+1}(\vecx)-\vecm^n(\vecx)}{\tau}\right)\right|\leq \tau \left|\partial_t \vecm_{h,\tau}(\vecx,t)\right|.\nonumber
\end{equation}
Therefore
\begin{equation*}
||\vecm_{h,\tau}-\vecm^-_{h,\tau}||_{L^2(\OT)}\leq \tau\left\Vert\partial_t \vecm_{h,\tau}\right\Vert_{L^2(\OT)}\rightarrow_{(h,\tau)\rightarrow 0} 0,
\end{equation*}
which entails that 
\begin{equation*}
\vecm^-_{h,\tau}\rightarrow_{(h,\tau)\rightarrow 0} \bar\vecm\text{ strongly in }L^2(\OT).
\end{equation*}
Moreover, on any tetrahedron $K$ of $\mathcal{T}_h$, and for any
$\vecu\in M_h$ one has, $\vecx_i^h$ being any vertex of $K$,
\begin{equation*}
\left||\vecu(\vecx)|-|\vecu(\vecx_i^h)|\right|^2\leq Ch^2|\nabla \vecu|^2,
\end{equation*}
(recall that $\nabla \vecu$ is constant on $K$), from which one deduces (since $|\vecm_{h,\tau}^-(\vecx_i^h)|=1$)
\begin{equation*}
  \int_{\OT}\left|1-|\vecm^-_{h,\tau}|\right|^2\dx\leq Ch^2||\nabla \vecm^-_{h,\tau}||^2_{L^2(\OT)}.
\end{equation*}
This shows that $|\bar\vecm(\vecx,t)|=1$ a.e. $(\vecx,t)\in \OT$.
\medskip

Eventually, from the fact that at each vertex $\forall i\in\{1,\cdots,N_h\}$
\begin{equation}
  |\vecm^{n+1}_i-\vecm^n_i-\tau \vecv^n_i|=|\vecm^n_i+\tau \vecv^n_i|-1\leq \half \tau^2|\vecv^n_i|^2,
\label{vvv}
\end{equation}
we derive
\begin{equation*}
  \left|\frac{\vecm^{n+1}_i-\vecm^n_i}{\tau}-\vecv^n_i\right|\leq \half \tau|\vecv_i^n|^2.
\end{equation*}
Appealing to \eqref{ineg} the latter entails that
\begin{equation}
  \left\Vert\partial_t \vecm_{h,\tau}-\vecv_{h,\tau}\right\Vert_{L^1(\OT)}\leq c^2 \tau ||\vecv_{h,\tau}||^2_{L^2(\OT)}\rightarrow_{(h,\tau)\rightarrow 0} 0.\nonumber
\end{equation}
This is sufficient to conclude that $\displaystyle \vecv=\partial_t
\bar \vecm$ in \eqref{21}.
\paragraph{General properties of interpolation operator}
Before we start with the penultimate step of proving convergence, let
us state some general properties of the nodal interpolation operator
which we repeatedly use in the sequel.  Up to dimension three, there
holds for any function $\varphi\in H^2(\Omega)\subset \mathcal
C^0(\bar \Omega)$
\begin{align}
  \label{JS:Interpol1}
  ||\varphi-\mathcal{I}_h(\varphi)||_{H^1(\Omega)} \leq C h ||\nabla^2
  \varphi||_{L^2\Omega}.
\end{align}
Since the basis functions are linear on each triangle one can deduce form \eqref{JS:Interpol1} that 
\begin{align}
  \label{JS:Interpol2}
  ||\vecm_{h,\tau}^-\times \tilde
  \tf-\mathcal{I}_h(\vecm_{h,\tau}^-\times\tilde\tf)||_{L^2([0,T],H^1)}
  \leq C h ||\vecm_{h,\tau}^-||_{H^1(\OT)}||\tf||_{W^{2,\infty}},
\end{align}
see \cite[p.7]{alougesAIMS}. 
\paragraph{Convergence to a solution of the LLG equation}
Having established the preliminary results above, we are now ready to
proceed with the proof of convergence: Test \eqref{JS:inter5} with
$\tf=\mathcal{I}_h( \vecm^-_{h,\tau}\times \tilde\tf)$ where
$\tilde\tf \in C_0^\infty(\OT)^3$. We recall that $\mathcal{I}_h$ is
the nodal interpolation, cf.\ \eqref{JS:interpol}. After suitable
integration in time 
we hence obtain from \eqref{JS:inter5.2} with the choice of
$\tf=\mathcal{I}_h( \vecm^-_{h,\tau}\times \tilde\tf)$ that
\begin{multline}
  \int_\OT \varphi_M(\vecm^-_{h,\tau}) \,\vecv_{h,\tau}\cdot
  \mathcal{I}_h(\vecm_{h,\tau}^-\times \tilde\tf) \dx\dt
  +\int_\OmT  \vecm_{h,\tau}^-\times \vecv_{h,\tau}\cdot \mathcal{I}_h(\vecm_{h,\tau}^-\times \tilde\tf)\dx\dt\\
  +\tfrac{\tau}{2} \int_\OmT (1+\rho(\tau))\,d^2\,\nabla
  \vecv_{h,\tau} \cdot \nabla \mathcal{I}_h(\vecm_{h,\tau}^-\times
  \tilde\tf) -\Hd(\vecv_{h,\tau})\cdot
  \mathcal{I}_h(\vecm_{h,\tau}^-\times \tilde\tf)\\
  -Q(\vece\cdot \vecv_{h,\tau})(\vece
  \cdot\mathcal{I}_h(\vecm_{h,\tau}^-\times \tilde\tf))
  \dx\dt\\
  =\int_\OmT -\,d^2\, \nabla \vecm_{h,\tau}^-\cdot \nabla
  \mathcal{I}_h(\vecm_{h,\tau}^-\times \tilde\tf)
  +\Hd(\vecm_{h,\tau})\cdot \mathcal{I}_h(\vecm_{h,\tau}^-\times
  \tilde\tf) \\
  +Q(\vece \cdot\vecm_{h,\tau})( \vece
  \cdot\mathcal{I}_h(\vecm_{h,\tau}^-\times \tilde\tf) ) +H_\ext \cdot
  \mathcal{I}_h(\vecm_{h,\tau}^-\times \tilde\tf)
  \dx\dt.\label{JS:inter7}
\end{multline}
Our goal is to pass to the limit $(\tau,h)\rightarrow 0$ in the latter
equation \eqref{JS:inter7} to recover the LLG equation
\eqref{disLLG}. As we shall see, the first and the third term on the
l.h.s.\ and the first term on the r.h.s.\ are a little bit subtle and have
to be treated with caution. The remaining contributions behave well
under the established convergence; this is particularly due to the
fact that $\Hd$ is $L^2$-continuous. For the second contribution on
the l.h.s. one further uses that the $L^\infty$ bound on $\vecm^-$
improves \eqref{20} to strong convergence in any $L^p$ with $1 <
p<+\infty$.  \medskip

Let's start with the first contribution on the l.h.s.
Observe that $|\varphi_M|$ is uniformly bounded. Moreover it holds
that $|\varphi_M-\alpha|\leq \frac{\tau M}{2}$. As long as $\tau
M\rightarrow 0 $ for $(h,\tau)\rightarrow 0$ the strong convergence of
$\vecm_{h,\tau}^-$ is sufficient to conclude that
\begin{equation} 
 \label{JS:con1}
  \int_\OT  \varphi_M(\vecm^-_{h,\tau}) \,\vecv_{h,\tau}\cdot \mathcal{I}_h(\vecm_{h,\tau}^-\times \tilde\tf) \dx\dt\rightarrow_{(h,\tau)\rightarrow 0}  \alpha \int_\OmT \vecv \cdot( \bar \vecm \times \tilde \tf)\dx\dt. 
\end{equation}
In fact, using the triangle inequality we find that
\begin{multline}
  \left|\int_\OT \varphi_M(\vecm^-_{h,\tau}) \,\vecv^n_{h,\tau}\cdot
    \mathcal{I}_h(\vecm_{h,\tau}^-\times \tilde\tf) \dx\dt
    -\alpha \int_\OmT \vecv \cdot( \bar \vecm \times \tilde \tf)\dx\dt\,\right|\\
  \leq \left|
    \int_\OT \varphi_M(\vecm^-_{h,\tau}) \,\vecv_{h,\tau}\cdot(
    \vecm_{h,\tau}^-\times \tilde\tf)\dx\dt-\alpha \int_\OmT \vecv
    \cdot( \bar \vecm \times \tilde\tf)\dx\dt\,\right|\\
  +\left|
    \int_\OT \varphi_M(\vecm^-_{h,\tau})
    \,\vecv_{h,\tau}\cdot((\vecm_{h,\tau}^-\times \tilde\tf)-
    \mathcal{I}_h(\vecm_{h,\tau}^-\times \tilde\tf))\dx\dt\,\right|. \label{triangle}
\end{multline}
The first term tends to zero since 
\begin{eqnarray*}
  & \varphi_M(\vecm^-)
  \rightarrow_{(h,\tau)\rightarrow 0} \alpha &\quad\text{ in }\quad L^\infty(\Om),\\
  & \vecm_{h,\tau}^-\rightarrow_{(h,\tau)\rightarrow 0} \bar \vecm
  \quad\text{ in }&\quad
  L^2(\OmT), \text{ and}\\
  & \vecv_{h,\tau}\rightarrow_{(h,\tau)\rightarrow 0}
  \vecv=\displaystyle\frac{\partial \bar\vecm}{\partial t}&\quad\text{ in }\quad L^2(\OmT)
  \end{eqnarray*}
as ${h,\tau}\rightarrow 0$.

Since $\varphi_M(\vecm^-_{h,\tau}) $ is uniformly bounded, we can evoke
\eqref{JS:Interpol2} to obtain that the second contribution tends to
zero. This establishes \eqref{JS:con1}.  \medskip
 
Let's turn to the next term in \eqref{JS:inter7}. Convergence in this
case essentially relies upon the estimate \eqref{bound1}. In fact,
appealing once again to \eqref{JS:Interpol2} we see that instead of
establishing
\begin{equation}
  \label{JS:10.1}
  \tfrac{\tau}{2} \,d^2\,\int_\OmT (1+\rho(\tau))\nabla \vecv_{h,\tau} \cdot \nabla \mathcal{I}_h(\vecm_{h,\tau}^-\times \tilde\tf)\dx
  \rightarrow_{(h,\tau)\rightarrow 0} 
  0,
\end{equation}
if suffices to establish
\begin{equation}
  \label{JS:10.2}
  \tfrac{\tau}{2} \,d^2\,\int_\OmT (1+\rho(\tau))\nabla \vecv_{h,\tau} \cdot \nabla (\vecm_{h,\tau}^-\times \tilde\tf)\dx
  \rightarrow_{(h,\tau)\rightarrow 0} 
  0,
\end{equation}
which follows obviously from \eqref{bound1} using Young's inequality.
\medskip

Finally, the convergence of the last term in \eqref{triangle} follows from the
orthogonality property of the cross product and \eqref{19}, \eqref{20}
by once again appealing to \eqref{JS:Interpol2} since
\begin{align}
  \tfrac{\tau}{2}\Big| \int_\OmT \nabla \vecm_{h,\tau}^-\cdot \nabla&
    \mathcal{I}_h(\vecm_{h,\tau}^-\times \tilde\tf) \dx\dt
    -\int_\OmT \nabla \bar \vecm\cdot \bar \vecm \times \nabla\tilde\tf\dx\dt\Big|\notag\\
  &\leq \tfrac{\tau}{2}\left| \int_\OmT \nabla \vecm_{h,\tau}^-\cdot
    \nabla \left(\mathcal{I}_h(\vecm_{h,\tau}^-\times \tilde\tf)
      -(\vecm_{h,\tau}^-\times \tilde\tf)\right) \dx\dt\right| \\&\qquad+
  \tfrac{\tau}{2}\left|\int_\OmT \nabla \vecm_{h,\tau}^-\cdot \nabla
    (\vecm_{h,\tau}^-\times \tilde\tf)
    -\nabla \bar \vecm\cdot \bar \vecm \times \nabla\tilde\tf\dx\dt\right|\notag\\
  &= \tfrac{\tau}{2}\left| \int_\OmT \nabla \vecm_{h,\tau}^-\cdot
    \nabla \left(\mathcal{I}_h(\vecm_{h,\tau}^-\times
      \tilde\tf)-(\vecm_{h,\tau}^-\times \tilde\tf)\right) \dx\dt\right|\\&\qquad
  + \tfrac{\tau}{2}\left|\int_\OmT \nabla \vecm_{h,\tau}^-\cdot(
    \vecm_{h,\tau}^-\times \nabla \tilde\tf)-\nabla \bar \vecm\cdot
    (\bar \vecm \times \nabla\tilde\tf)\dx\dt\right|. \label{JS:11}
\end{align}
\paragraph{Energy estimate}
We finally establish the energy estimate.  
From \eqref{estimenerg} we deduce that $\forall n\in\{0,\cdots,N_h\}$
\begin{align}
  \mathcal{E}(\vecm^{n+1})-\mathcal{E}(\vecm^n)&\leq-2\alpha
  \tau \int_{\Omega} \varphi(\vecm^n)\left|\vecv^n\right|^2\dx+2\tau\int_\Omega \bH\left(\vecm^n\right)\cdot \vecv^n\dx\notag\\
  &\quad+\tau^2\int_\Omega
  \frac{\partial\bH}{\partial_\vecm}(\vecv^n)\cdot \vecv^n\dx
  -\tau^2\rho(\tau)\,d^2\,\int_\Omega\left|\nabla \vecv^n\right|^2\dx\notag\\
  &\quad-\int_\Omega
  \bH\left(\vecm^{n+1}+\vecm^n\right)\cdot(\vecm^{n+1}-\vecm^n)\dx,\label{JS:energ}
\end{align}
cf. \eqref{E(m)}. Let us introduce another short-hand notation for the remaining effective field, namely $\bH^n=\bH\left(\vecm^n\right)$.
We consider the contributions in \eqref{JS:energ} separately and start with the observation that
\begin{align*}
  2\tau\int_\Omega \bH^n\cdot \vecv^n\dx&-\int_\Omega (\bH^{n+1}+\bH^{n})\cdot(\vecm^{n+1}-\vecm^n)\dx\\
  &= 2\int_\Omega \bH^n\cdot (\vecm^{n+1}-\vecm^n-\tau
  \vecv^n)\dx+\int_\Omega (\bH^{n+1}-\bH^{n})\cdot(\vecm^{n+1}-\vecm^n)\dx.
\end{align*}
Hence due to  \eqref{ineg2} and \eqref{vvv} combined with \eqref{ineg}
\begin{align}
  \left|2\tau\int_\Omega \bH^n\cdot \vecv^n\dx-\int_\Omega
    (\bH^{n+1}+\bH^{n})\cdot(\vecm^{n+1}-\vecm^n)\dx\right|\dx\leq C
  \tau^2(||\vecv^n||_{L^2}||\vecv^n||_{L^4}+||\vecv^n||^2_{L^2}) \label{energ6}
\end{align}
In order to bound the stray-field contribution we have employed \eqref{estimcst} with $p=4$. 
The contributions in the second line of the r.h.s. of \eqref{JS:energ} are of higher order in $\tau$. The first term can be easily bounded using Young's inequality:
\begin{align}
  \left| \int_\Omega \frac{\partial\bH}{\partial \vecm}(\vecv^n)\cdot
    \vecv^n\dx\right|\leq C ||\vecv^n||^2_{L^2}. \label{energ7}
\end{align}
Plugging in \eqref{energ6} and \eqref{energ7} into \eqref{JS:energ} yields that
\begin{align*}
  \mathcal{E}(\vecm^{n+1})-\mathcal{E}(\vecm^n)&+2\tau\int_{\Omega} \varphi(\vecm^n)\left|\vecv^n\right|^2\dx\\
  &\leq C \tau^2(||\vecv^n||_{L^2}+||\vecv^n||^2_{L^2}+||\vecv^n||_{L^2}||\vecv^n||_{L^4}+\rho(\tau)\,d^2\,||\nabla \vecv^n||^2_{L^2})\\
  &\leq
  C'\tau^2(||\vecv^n||_{L^2}+||\vecv^n||^2_{L^2}+||\vecv^n||_{L^2}||\nabla
  \vecv^n||_{L^2}+\rho(\tau)||\nabla \vecv^n||^2_{L^2}),
\end{align*}
where $C$ denotes a generic constant. Here we made use of the classical Sobolev
embedding
$$
||\vecv^n||_{L^4}\leq C||\nabla \vecv^n||_{L^2}\,.
$$
Summing from $n=0$ to $N-1$ leads to
\begin{multline*}
  \mathcal{E}(\vecm(N\tau))-\mathcal{E}(\vecm(0))+
  \int_\OT \varphi_M(\vecm_{h,\tau}^-)|\vecv_{h,\tau}|^2\dx\dt\\
  \leq C \tau
  (||\vecv_{h,\tau}||_{L^2}+||\vecv_{h,\tau}||^2_{L^2}+||\vecv_{h,\tau}||_{L^2}||\nabla
  \vecv_{h,\tau}||_{L^2}+\rho(\tau)||\nabla \vecv_{h,\tau}||^2_{L^2}).
\end{multline*}
We are now ready to pass to the limit. Noticing once again that ${\tau}||\nabla \vecv^n||_{L^2(\OT)})$ is uniformly bounded from \eqref{bound1} we derive that
\begin{equation}
\mathcal{E}(\vecm(N\tau))-\mathcal{E}(\vecm(0))  +\alpha\int_0^T\int_{\Omega}\left|\vecv^n\right|^2\dx\,\dt \leq 0.
\end{equation}


\bibliographystyle{elsarticle-num} 







\end{document}